\documentclass{amsart}
\usepackage{graphicx}
\usepackage{amsmath}
\usepackage{amsthm}
\usepackage[dvips]{color}
\usepackage{multicol}
\usepackage{amsfonts}
\usepackage{amssymb, amscd}
\newtheorem{thm}{Theorem}[section]
\newtheorem{cor}[thm]{Corollary}
\newtheorem{lemma}[thm]{Lemma}
\newtheorem{prop}[thm]{Proposition}
\newtheorem{example}[thm]{Example}
\theoremstyle{definition}
\newtheorem{defn}[thm]{Definition}
\theoremstyle{remark}
\newtheorem{rem}[thm]{Remark}
\numberwithin{equation}{section}

\newcommand{\K}{\mathbb K}

\newcommand{\A}{\mathcal{A}}

\begin{document}

\title[Bialgebra structures of 2-associative algebras]{Bialgebra structures of 2-associative algebras}%

\author{K. DEKKAR and A. MAKHLOUF}%
\address{Khadra DEKKAR, Universit\'{e}
de S\'{e}tif, Facult\'e des Sciences et Techniques,  Algeria}
\email{dekkar{\_}na@yahoo.fr }%
\address{Abdenacer MAKHLOUF {(\tiny{corresponding author})}, Universit\'{e} de Haute Alsace, Laboratoire de Math\'{e}matiques, Informatique et
Applications, Mulhouse, France} \email{Abdenacer.Makhlouf@uha.fr}


\begin{abstract}

This work is devoted to  study  new bialgebra structures related to
2-associative algebras. A 2-associative algebra  is a vector space
equipped with two associative
    multiplications. We discuss the notions of 2-associative bialgebras,
    2-bialgebras and 2-2-bialgebras. The first structure was revealed by J.-L. Loday and
M. Ronco in an analogue of a Cartier-Milnor-Moore theorem,
 the second was suggested by Loday and the third is a variation of
 the second one.
 The main results of this paper are the construction of   2-associative
 bialgebras,  2-bialgebras and 2-2-bialgebras starting from an
  associative algebra and  the classification of these structures in
  low dimensions.

\end{abstract}
\maketitle
\section*{Introduction}
The aim of this work is to  study some new  algebraic structures
related to  bialgebras  and 2-associative algebras. The motivation
comes from an extension to non-cocommutative situation of the
Cartier--Milnor--Moore theorem. Namely, we discuss 2-associative
bialgebras, 2-bialgebras and 2-2-bialgebras. We provide, using
Kaplansky's bialgebras (see \cite{Ka73}), a construction of
  $2$-associative bialgebras and  2-bialgebras, starting from any two
  $n$-dimensional
associative algebras.  Also, we establish some properties and
classifications in low dimensions.

A \textit{bialgebra}  is a $\K$-vector space $V$, where $\K$ is a
field, equipped with an algebra structure given by a multiplication
$\mu$ and a unit $\eta$ and a
 coalgebra structure given by a comultiplication
$\Delta$ and a counit $\varepsilon$, such that there is a
compatibility condition between the two structures expressed by the
fact that $\Delta$ and $\varepsilon$ are algebra morphisms, that is
for $x,y\in V$
$$\Delta(\mu(x\otimes y))=\Delta(x)\bullet \Delta(y)\quad \text{and}
\quad \varepsilon(\mu(x,y))=\varepsilon(x)\varepsilon (y).$$
The multiplication $\bullet $ on $V\otimes V$  is the usual
multiplication on tensor product $$ \left( x\otimes
y\right) \bullet \left( x^{\prime }\otimes y^{\prime }\right) =\mu
\left(x\otimes x^{\prime } \right) \otimes \mu \left( y\otimes y^{\prime
}\right)
$$
 We assume also that the unit is sent to the unit by the
comultiplication. A bialgebra is said to be a \textit{Hopf algebra}
if the identity on $V$ has an inverse for the convolution product
defined by
\begin{equation}\label{conv} f\star g := \mu
\circ (f\otimes g)\circ \Delta\
\end{equation}

 Let  $(C,\Delta,\varepsilon )$ be a graded coalgebra over a field
$\K$, that is $\mathcal{C}=\oplus_{k>0}\mathcal{C}^k$  such that
$\Delta(\mathcal{C}^k)\subseteq  \sum_{ i+j=k} \mathcal{C}^i \otimes
\mathcal{C}^j$ and $\varepsilon(C^k) = 0 \ \forall k \neq 0.$  The
coalgebra is said to be \textit{graded connected} if in addition
$\mathcal{C}^0 \cong \K.$ A graded coalgebra
$\mathcal{C}=\oplus_{k>0}\mathcal{C}^k$ is said to be
\textit{cofree} if it satisfies the following universal property.
Given a graded coalgebra $\mathcal{Q}=\oplus_{k>0}\mathcal{Q}^k$ and
a linear map $\varphi : \mathcal{Q} \rightarrow \mathcal{C}^1$ with
$\varphi(\mathcal{Q}^k) = 0$ when $k \neq 1$, there is a unique
morphism of graded coalgebras $\widehat{\varphi} : \mathcal{Q}
\rightarrow \mathcal{C}$ such that $\pi \circ
\widehat{\varphi}=\varphi$ where $\pi : \mathcal{C} \rightarrow
\mathcal{C}^1$ is the canonical projection. We refer to (see
\cite{Kas95}, \cite{Mo92}, \cite{Ma05}) for the basic knowledge
about Hopf algebras.

  The \textit{Cartier--Milnor--Moore theorem} (see \cite{MM},\cite{Ca})
  states that, over a field of characteristic
  zero, any connected cocommutative bialgebra $\mathcal{H}$ is
  of the form $U(Prim(\mathcal{H}))$, where the primitive part
  $Prim(\mathcal{H})$ is viewed as a Lie algebra, and $U$ is the
  universal enveloping  functor. Thus there is an equivalence of categories  between the
  cofree cocommutative  bialgebras and  Lie algebras.

  In (see \cite{Lo-Ro04}), J.-L.
  Loday and M. Ronco change the category of connected cocommutative bialgebras
  for the category of 2-associative
  bialgebras, the Lie algebras by a
  $B_\infty$-algebras and the universal enveloping functor by a
  functor $U2$ $$U2 : B_\infty \text{-algebras} \longrightarrow \text{2-associative
  algebras.}$$ A non-differential
    $B_\infty$-algebra is defined by $(p+q)$-ary operations for any
    pair of positive integers $p,q$ satisfying some relations. It may be viewed also as a deformation of
    the shuffle algebra $T^{sh}(V)$ where $V$ is the underlying vector space.
      Loday and Ronco proved that if $\mathcal{H}$ is a bialgebra over a field $\K$
the following assertions are equivalent
\begin{enumerate}
\item $\mathcal{H}$ is a connected $2$-associative bialgebra,
\item $\mathcal{H}$ is isomorphic to $U2(Prim\ \mathcal{H})$ as a
$2$-associative bialgebra,
\item $\mathcal{H}$ is cofree among the connected coalgebras.
\end{enumerate}
Therefore, any cofree Hopf algebra is of the form $U2(R)$ where $R$
is $B_\infty \text{-algebras}$. J.-L. Loday gives in (see
\cite{Loday-GB}) a general framework  to work with triples of
operads including the previous case.

The paper is organized as follows: Section 1 is dedicated to
introduce the definitions  of 2-associative algebras, 2-associative
bialgebras, infinitesimal bialgebras, 2-bialgebras and
2-2-bialgebras, and to describe some basic properties of these
objects. In Section 2, we recall how I. Kaplansky (see \cite{Ka73})
constructed from an associative algebra two  bialgebras. We show
that they leads to a large class of 2-associative bialgebras,
2-bialgebras and 2-2-bialgebras. In Section 3, we establish
classifications in dimensions 2 and 3 of bialgebra, bialgebras which
carry also infinitesimal bialgebra structure. Therefore, we
enumerate and describe 2 and 3-dimensional 2-associative bialgebras,
2-bialgebra and 2-2-bialgebras.

\section{Definitions and Properties}
In this section, we introduce the definitions of  bialgebras,
2-bialgebras and 2-2-bialgebras, which requires the notions of
$2$-associative algebra, bialgebra and infinitesimal bialgebra.
\subsection{$2$-Associative algebras}
We review first briefly the 2-associative algebra structure and some
properties, in particular concerning the free 2-associative algebra.
A theory of this structure was developed by J.-L. Loday and M. Ronco
in (see \cite{Lo-Ro04}).
\begin{defn}
A \emph{2-associative algebra} over $\mathbb{K}$ is a vector space equipped with two associative operations. A 2-associative algebra is said to be unital if there is a unit $1$ which is a unit for both operations.
\end{defn}

One may define the notions of 2-associative monoid, 2-associative
group, 2-associative  monoidal category, 2-associative operad ...
The free 2-associative algebra over the vector space $V$ is the
2-associative algebra $2as(V)$ such that any linear map from $V$ to
a 2-associative algebra $\mathcal{A}$ has a unique extension to  a
homomorphism of  2-associative algebras from $2as(V)$ to
$\mathcal{A}$. An explicit description in terms of planar trees is
given in (see \cite{Lo-Ro04}). The free 2-associative algebra on one
generator can be identified with the non-commutative polynomials
over the planar (rooted) trees. The operad of 2-associative algebras
is studied in (see \cite{Lo-Ro04}), and it turns out that it is a
Koszul operad.

\subsection{Infinitesimal bialgebras}
An infinitesimal bialgebra is a bialgebra where the compatibility
condition is modified. The comultiplication is no more an algebra morphism. The condition is
$
\Delta \circ \mu = (\mu \otimes id)\circ(id \otimes \Delta)+(id
\otimes \mu)\circ( \Delta \otimes id).
$
 This structure was introduced first by
S.~Joni and G.-C.~Rota in (see \cite{JR}).  The basic theory was
developed by M.~Aguiar in (see \cite{Aguiar2000} \cite{Aguiar2001}).
He also showed their intimate link to Rota--Baxter algebras, Loday's
dendriform algebras, pre-Lie structure and introduced the
associative classical Yang--Baxter equation (see \cite{Aguiar}). The
following definition was introduced by J.-L.~Loday and M.~Ronco.

 \begin{defn}A \emph{unital infinitesimal bialgebra} $%
\left( V,\mu ,\eta,\Delta,\varepsilon \right) $ is a vector space
$V$ equipped with a
unital associative multiplication $\mu $ and a counital coassociative comultiplication $%
\Delta $ which are related by the unital infinitesimal relation
\begin{equation}\label{inf}\Delta \left( \mu (x\otimes y) \right) =\left( x\otimes 1\right) \bullet
\Delta \left( y\right) +\Delta \left( x\right) \bullet \left(
1\otimes y\right) -x\otimes y\end{equation} where $1=\eta(1)$ and
$x,y\in V.$
\end{defn}

\begin{rem}
The unital infinitesimal relation \ref{inf} may be written as follows :
\begin{equation}\label{inf2}
\Delta \circ \mu = (\mu \otimes id_V)\circ(id_V \otimes \Delta)+(id_V
\otimes \mu)\circ( \Delta \otimes id_V )-id_V \otimes id_V
\end{equation}
or
$$\Delta \circ \mu =\left( \mu \otimes \mu \right)\circ \left(
id_{V}\otimes \tau\otimes id_{V}\right) \mathit{\circ }\left(
i_{1}\otimes \Delta +\Delta \otimes i_{2}-i_{1}\otimes
i_{2}\right)$$ where $\tau$ is the usual flip and $$ i_{1}:  V
\rightarrow V\otimes V \ \
 \ x  \rightarrow  x\otimes 1$$
$$
i_{2}:  V  \rightarrow  V\otimes V \ \
 y   \rightarrow 1\otimes y$$
The maps  $i_{1}$ and $i_{2}$ are algebra morphisms.

\end{rem}

\begin{rem}
The condition \ref{inf2} applied to $1\otimes1$ implies that $\Delta
(1)=1 \otimes 1.$
\end{rem}
The polynomial algebra and the tensor algebra may be endowed with
infinitesimal bialgebra structure, (see \cite{Lo-Ro04}).
\begin{example}
The polynomial algebra $\K [x]$ is a unital infinitesimal bialgebra
with    a comultiplication $\Delta$ defined by $\Delta(x^n)= \sum
_{p=0}^{n} x^p \otimes x^{n-p}$.
\end{example}

\begin{example}The tensor algebra over a vector space $V$ which is the
space $$T(V) = \K \oplus V \oplus  V^{\otimes 2} \oplus \cdots
\oplus V^{\otimes n} \oplus \cdots $$ equipped with  the
concatenation multiplication given by $v_1\cdots v_i \otimes v_{i+1}
\cdots v_n \mapsto v_1 \cdots v_i v_{i+1} \cdots v_n \ , $ and the
deconcatenation comultiplication given by  $ \Delta( v_1 \cdots v_n
) = \sum_{i=0}^{i=n}v_1\cdots v_i \otimes v_{i+1} \cdots v_n \ $ is
a unital infinitesimal bialgebra.

In (see \cite{Lo-Ro04}), the authors showed that any connected
unital infinitesimal bialgebra is isomorphic to $T(V)$ for some
space $V$.
\end{example}

We may extend the  unital  infinitesimal bialgebra structure to a situation
where the multiplication
$\mu$ is no more unital and the comultiplication $\Delta$ is no more
counital by considering a
triple $ \left( V,\mu ,\Delta \right) $  with the infinitesimal
compatibility condition
(\ref{inf2}).

The unital infinitesimal relation differs from the infinitesimal
relation used by S.~Joni and G.-C.~Rota  by the presence of the term
$-x\otimes y$. The original definition assumes that the
comultiplication is a  derivation. One may generalize the unital
infinitesimal condition, with $\theta\in \K$  (see \cite{KEF}) to
\begin{equation}\label{infG}
\Delta \circ \mu = (\mu \otimes id_V)(id_V \otimes \Delta)+(id_V
\otimes \mu)( \Delta \otimes id_V )-\theta\  id_V \otimes id_V.
\end{equation}
The case $\theta =0$ is the relation defined by Joni and Rota, and
$\theta =1$ recovers the notion of unital infinitesimal bialgebra.

 The convolution product (\ref{conv}) defined by a unital
infinitesimal bialgebra or a generalized infinitesimal algebra
(condition \ref{infG}) still endows the set $Hom(V,V)$ of
$\K$-linear endomorphisms of $V$, with a structure of associative
algebra where $\varepsilon \circ \eta $ is the unit. Any
infinitesimal bialgebra induces over the algebra $(Hom(V,V),\circ
)$, where $\circ$ denotes the composition of maps, a Rota--Baxter
structure (see \cite{KEF}), that is a $\K$-linear map $\phi :
Hom(V,V)\rightarrow Hom(V,V)$ satisfying, for all $f,g\in Hom(V,V)$,
$$\phi (f) \circ \phi(g)-\phi(f\circ g)=\phi(\phi(f)\circ g+f\circ\phi(g)).$$
One may consider $\phi (f)=id_V \star f$ or  $\phi (f)= f\star id_V
$, where $\star$ denotes the convolution product.

Aguiar  showed also that to any infinitesimal bialgebra $(\theta
=0)$, corresponds a preLie algebra (see \cite{Aguiar2000}[theorem
3.2]). A preLie algebra  is a vector space $V$ together with a
bilinear map $m$ satisfying, for all $x,y,z\in V$, the following
condition
$$m(x, m(y,z))-m(m(x, y),z)=m(y, m(x,z))-m(m(y, x),z).
$$
The commutator $[x,y]=m(x,y)-m(y,x)$ defines a Lie algebra on $V$.
The infinitesimal algebra $(V,\mu,\Delta)$ induces the preLie multiplication defined, using Sweedler notation, for $x,y\in V$ by
$$
m(x,y)=\sum_{(y)}{\mu(\mu(y_{(1)}\otimes x)\otimes y_{(2)})}.
$$

 \subsection{ $2$-Associative Bialgebras}
A 2-associative bialgebra is defined as follows.

 \begin{defn}
 A \emph{2-associative bialgebra}  $\mathcal{B}2as=\left( V,\mu_1
,\mu_2,\eta,\Delta,\varepsilon \right) $ is a vector space $V$
equipped with two multiplications $\mu_1 $ and $\mu_2$, a unit
$\eta$, a comultiplication $\Delta $ and a counit $\varepsilon$ such
that
\begin{enumerate}
\item $ \left( V,\mu_1, \eta,\Delta ,\varepsilon\right)$ is a bialgebra,

\item $ \left( V,\mu_2 ,\eta,\Delta ,\varepsilon\right)
$ is a unital infinitesimal bialgebra.
\end{enumerate}
\end{defn}
One may define in a similar way a  2-associative Hopf algebra and a
unital infinitesimal Hopf algebra.

Let $\left( V,\mu_1 ,\mu_2%
 ,\eta,\Delta ,\varepsilon \right) $ and $\left(
V^{\prime },\mu_1 ^{\prime },\mu^{\prime
},\eta_2^{\prime },\Delta ^{\prime },\varepsilon ^{\prime }\right) $ be
two 2-associative bialgebras. A linear map $f: V\rightarrow
V^{\prime }$ is a morphism of 2-associative bialgebras if
$$
\mu_1^{\prime} \circ \left( f\otimes f\right) =f\circ \mu_1, \quad  \quad\mu_2^{\prime }\circ \left( f\otimes f\right) =f\circ
\mu_2,
 \quad  \quad f\circ \eta  =\eta^{\prime },
\quad $$ $$
(f\otimes f)\circ \Delta =\Delta^{\prime } \circ f\quad \text {and} \quad \varepsilon^{\prime } =\varepsilon \circ f$$

 \subsection{$2$-Bialgebras and $2$-$2$-Bialgebras}
In this section, we set up the definitions of $2$-bialgebra and
$2$-$2$-bialgebra, and we give some basic properties and examples.

 \begin{defn}
 A \emph{2-bialgebra} $\mathcal{B}2=\left( V,\mu_1 ,\mu _2,\eta ,\Delta _1,\Delta _2 ,\varepsilon _1 ,\varepsilon _2\right) $ is a
vector space V equipped with two multiplications $\mu_1 $ , $\mu
_2$, one unit $\eta $,  two comultiplications $\Delta_1 ,$
$\Delta_2,$  two counits $\varepsilon_1 ,\varepsilon _2
$  such that
 $ \left( V,\mu_1 ,\eta ,\Delta_1 ,\varepsilon_1 \right) $, $\left( V,\mu_2,\eta ,\Delta _2,\varepsilon _2\right) $, $ \left( V,\mu_1 ,\eta ,\Delta_2
,\varepsilon_2 \right) $ and $\left( V,\mu_2,\eta
,\Delta_1 ,\varepsilon_1 \right) $ are bialgebras.
\end{defn}

\begin{rem}
The  condition could be expressed by
\begin{enumerate}
\item $ \mu_1 $ is  compatible with $\Delta_1$ and $\Delta_2.$

\item $ \mu _2$ is  compatible with $\Delta _2$ and  $\Delta_1 .$

\item $\mathit{\varepsilon_1 \circ \mu _1=\varepsilon_1 \otimes
\varepsilon_1 }$  and $ \mathit{\varepsilon_1 \circ \mu _2=\varepsilon_1 \otimes
\varepsilon_1 }$

\item $\mathit{\varepsilon }_2\mathit{\circ \mu_1 =\varepsilon_2 }%
\otimes\mathit{ \varepsilon }_2$ and $ \mathit{\varepsilon }_2\mathit{\circ \mu_2 =\varepsilon_2 }%
\otimes\mathit{ \varepsilon }_2$
\end{enumerate}
\end{rem}

Note that there is no relation assumed between $\Delta_1$ and
$\varepsilon_2$  (resp. $\Delta_2$ and $\varepsilon_1$).

 A 2-bialgebra is called of \emph{type} (1-1) (resp. of type
(2-2)) if the two multiplications and the two comultiplications are
identical (resp. distinct).

A 2-bialgebra is called of \emph{type} (1-2) (resp. of type (2-1))
if the two multiplications are identical (resp. distinct) and the
two comultiplications are distinct (resp. identical).

Let  $\left( V,\mu_1 ,\mu_2,\eta,\Delta_1 ,\Delta_2,\varepsilon_1,\varepsilon_2\right) $ and
$\left( V^{\prime },\mu_1^{\prime } ,\mu_2^{\prime },\eta^{\prime },\Delta_1^{\prime } ,\Delta_2^{\prime },\varepsilon_1^{\prime },\varepsilon_2^{\prime }\right) $
 be two 2-bialgebras. A linear map $f: V \rightarrow
V^{\prime}$ is a 2-bialgebra morphism  if
$$\mu_1 ^{\prime }\circ \left( f\otimes f\right) =f\circ \mu_1,
\quad \quad
\mu_2^{\prime }\circ \left( f\otimes f\right) =f\circ
\mu_2, \quad \quad  f\circ \eta =\eta ^{\prime },$$ $$
(f\otimes f)\circ \Delta_1 =\Delta_1 ^{\prime }\circ f, \quad %
\varepsilon_1 =\varepsilon_1 ^{\prime }\circ f, \quad
(f\otimes f)\circ \Delta_2=\Delta_2^{\prime }\circ f,
\quad  \varepsilon_2=\varepsilon_2^{\prime
}\circ f.
$$

\begin{example}
  Let $\left( V,\mu _{1},\eta,\Delta ,\varepsilon \right) $ and $\left( V,\mu _{2},,\eta,\Delta ,\varepsilon \right) $ be two bialgebras
over the vector space spanned by $\left\{ 1,x,y\right\} $, with $,\eta (1)=1$, and
defined by

$$
 \mu _{1}\left( x\otimes x\right) =x,\quad \mu _{1}\left( y\otimes
y\right) =y,\quad \mu _{1}\left( x\otimes y\right) =\mu _{1}\left(
y\otimes x\right) =0, $$
 $$\mu _{2}\left( x\otimes x\right) =x, \quad \mu _{2}\left( y\otimes
y\right) =\mu _{2}\left( x\otimes y\right) =\mu _{2}\left( y\otimes
x\right) =0,  $$
 $$\Delta \left( 1\right) =1\otimes 1,\quad \Delta \left( x\right)
=x\otimes x,\quad \Delta \left( y\right) =y\otimes 1+1\otimes y,$$
$$\varepsilon \left( 1\right) =1,\quad \varepsilon \left( x\right) =1,\quad \varepsilon \left(
y\right) =0.$$

We assume that $\eta$, , $\eta (1)=1$,  is a unit for both $\mu_1$ and $\mu_2$.

Then $\left( V,\mu _{1},\mu _{2},\eta, \Delta ,\Delta ,\varepsilon \right)
$ is a 2-bialgebra.

\begin{rem}  Let $\left( V,\mu,\eta ,\Delta ,\varepsilon \right) $ be a bialgebra then $$
\left( V,\mu ,\mu,\eta ,\Delta ,\Delta ,\varepsilon \right) \quad \text{ and }\quad \left(
V,\mu ,\mu ^{op},\eta,\Delta ,\Delta ^{cop},\varepsilon \right) $$ are
2-bialgebras, where $\mu ^{op}(x\otimes y)=\mu (y\otimes x)$ and $\Delta ^{cop}(x)=\tau\circ \Delta (x)$, with $\tau (x\otimes y)=y\otimes x.$   The first 2-bialgebra is of type (1,1) and the second
is of type (2,2) if the bialgebra is neither commutative nor
cocommutative.
\end{rem}

\end{example}

We introduce here the notion of 2-2-bialgebra, which is a variation of the previous structure.

 \begin{defn}
 A \emph{2-2-bialgebra} $\mathcal{B}22=\left( V,\mu_1 ,\mu _2,\eta, \Delta_1
,\Delta_2,\varepsilon_1 ,\varepsilon_2 \right) $ is
a vector space V equipped with two multiplications $\mu_1 $, $\mu
_2$, two comultiplications $\Delta_1 ,$ $\Delta_2$ ,
two counits $\varepsilon_1 ,\varepsilon_2$ and one unit $\eta
$, such that

\begin{enumerate}

\item $ \left( V,\mu_1 ,\eta ,\Delta_1 ,\varepsilon_1 \right) $\textit{\ and }$%
\left( V,\mu_2,\eta ,\Delta_2,\varepsilon_2 \right) $ are  bialgebras,

\item $ \left( V,\mu_1 ,\eta , \Delta_2 ,\varepsilon_2 \right) $\textit{\ and }$%
\left( V,\mu_2,\eta , \Delta_1 ,\varepsilon_1 \right) $ are unital
infinitesimal bialgebras.

\end{enumerate}
\end{defn}

 A 2-bialgebra is called of \textit{type} (1-1) (resp. of \textit{type}
(2-2)) if the two multiplications and the two comultiplications are
identical (resp. distinct).

A 2-bialgebra is called of \textit{type} (1-2) (resp. of
\textit{type} (2-1)) if the two multiplications are identical (resp.
distinct) and the two comultiplications are distinct (resp.
identical).

The definition of  2-2-bialgebra morphism
  is similar to 2-bialgebra morphism.

\begin{example}
We provide the following 3-dimensional  example of 2-2-bialgebra of
type (2,2). We consider a basis $\{e_1,e_2,e_3\}$  of
$\mathbb{K}^3$. The 7-uple $\left( \mathbb{K}^3,\mu _{1},\mu
_{2},\eta,\Delta _{1},\Delta _{2},\varepsilon _{1},\varepsilon
_{2}\right), $ is defined as follows :

the unit for both multiplications $\eta(1)=e_1$,

the multiplications
\begin{align*}
\mu _{1}\left( e_{1}\otimes e_{i}\right) &=\mu _{1}\left(
e_{i}\otimes e_{1}\right) =e_{i} \ \ i=1,2,3,\\
\mu _{1}\left( e_{j}\otimes e_{2}\right) &=\mu _{1}\left(
e_{2}\otimes e_{j}\right) =e_{j}\ \ j=2,3,\\
\mu _{1}\left(
e_{3}\otimes e_{3}\right) &=e_{3}.
\end{align*}
\begin{align*}
 \mu _{2}\left( e_{1}\otimes e_{i}\right) &=\mu _{2}\left(
e_{i}\otimes e_{1}\right) =e_{i} \ \ i=1,2,3,\\
\mu _{2}\left( e_{j}\otimes e_{2}\right) &=\mu _{2}\left(
e_{2}\otimes e_{j}\right) =e_{j}\ \ j=2,3,~\\
\mu _{2}\left(e_{3}\otimes e_{3}\right) &=0,
\end{align*}
the comultiplications and counits
\begin{align*}
 \Delta _{1}\left( e_{1}\right) &=e_{1}\otimes
e_{1};\\
 \Delta _{1}\left( e_{2}\right) &=e_{1}\otimes
e_{2}+e_{2}\otimes e_{1}-e_{2}\otimes e_{2};\\
 \Delta _{1}\left(e_{3}\right)& =e_{1}\otimes e_{3}+e_{3}\otimes e_{1}-e_{2}\otimes
e_{3};\
\end{align*}
$$\varepsilon _{1}\left( e_{1}\right) =1;\ \varepsilon
_{1}\left( e_{2}\right) =0;\ \varepsilon _{1}\left(
e_{3}\right) =0,$$
\begin{align*}
 \Delta _{2}\left( e_{1}\right)& =e_{1}\otimes
e_{1};\\
 \Delta _{2}\left( e_{2}\right) &=e_{1}\otimes
e_{2}+e_{2}\otimes e_{1}-e_{2}\otimes e_{2};\\
 \Delta _{2}\left(
e_{3}\right) &=e_{1}\otimes e_{3}+e_{3}\otimes e_{1}-e_{3}\otimes
e_{2};\
\end{align*}
$$\varepsilon _{2}\left( e_{1}\right) =1;\ \varepsilon
_{2}\left( e_{2}\right) =0;\ \varepsilon _{2}\left(
e_{3}\right) =0.$$
\end{example}

\section{Constructions} In this Section, we use Kaplansky's
constructions of bialgebras (see \cite{Ka73}) to built 2-associative
bialgebras, 2-bialgebras and 2-2-bialgebras. First we recall briefly
the definitions of these  bialgebras.
\begin{prop}\label{kap1}
 Let $\mathcal{A}=(V,\mu,\eta )$ be a unital associative algebra (where   $e_{2}:=\eta (1)$
 being the unit). Let $\widetilde{V}$ be the vector space spanned by $V$ and $e_{1}$, $\widetilde{V}=span(V,e_{1})$.

 We have the bialgebra $\mathcal{K}_{1}(\mathcal{A}):=(\widetilde{V},\mu_1,\eta_1,\Delta_1,\varepsilon_1 )$ where

the multiplication  $\mu_{1}$ is  defined by :
$$\mu_{1}\left( e_{1}\otimes x\right) =\mu_{1} \left(
x\otimes e_{1}\right) =x \ \ \ \ \ \forall x\in \widetilde{V}$$
$$\mu_{1} \left( x\otimes y\right) =\mu  \left( x\otimes y\right)
\ \ \ \ \ \ \ \ \ \ \  \forall x,y\in V, $$

the unit $\eta_1$ is given by $\eta_1 (1)=e_1$,

the comultiplication $\Delta _{1}$ is defined by :
$$\Delta _{1}\left( e_{1}\right) =e_{1}\otimes
e_{1}$$
$$\Delta _{1}\left( x\right) =x\otimes e_{1}+e_{1}\otimes
x-e_{2}\otimes x \ \ \ \ \ \ \ \forall x\in V $$

and the counit $\varepsilon_1 $ is defined by : $$\ \varepsilon_1
\left( e_{1}\right) =1 \textit{, } \varepsilon_1 \left( x\right) =0
\ \ \ \ \ \forall x\in V.$$
\end{prop}
\begin{proof}
Straightforward, (see \cite{Ka73}).
\end{proof}
\begin{rem}
The previous construction may be done even with a nonunital algebra.
\end{rem}
The second type bialgebra constructed by Kaplansky is given by the following proposition.
\begin{prop}\label{kap2}
Let $\mathcal{A}=(V,\mu,\eta )$ be a unital associative algebra (where   $e_{2}:=\eta (1)$
 being the unit). Let $\widetilde{V}$ be the vector space spanned by $V$ and $e_{1}$, $\widetilde{V}=span(V,e_{1})$.

 We have the bialgebra $\mathcal{K}_{2}(\mathcal{A}):=(\widetilde{V},\mu_2,\eta_2,\Delta_2,\varepsilon_2 )$ where

 the multiplication $\mu_{2}$ is
defined by:
\begin{eqnarray*}\mu_{2}\left( e_{1}\otimes x\right) &=&\mu_{2} \left(
x\otimes e_{1}\right) = x  \quad \forall x\in \widetilde{V}\\
\mu_{2} \left( x\otimes y\right) &=&\mu  \left( x\otimes
y\right)\quad
  \forall x,y\in V \end{eqnarray*}

  the unit $\eta_2$ is given by $\eta_2 (1)=e_1$,

  the comultiplication $\Delta
_{2}$ defined by :
\begin{eqnarray*}
\Delta _{2}\left( e_{1}\right) &=& e_{1}\otimes e_{1},\\
 \Delta_{2}\left( e_{2}\right) &=& e_{2}\otimes e_{1}+e_{1}\otimes
e_{2}-e_{2}\otimes e_{2}\\
\Delta _{2}\left( x\right) &=&\left( e_{1}-e_{2}\right) \otimes
x+x\otimes \left( e_{1}-e_{2}\right) \quad \forall x\in V\setminus
\{e_2\}, \end{eqnarray*}
and the counit defined by $\varepsilon $:
$$\ \varepsilon_2 \left( e_{1}\right) =1 \textit{, } \varepsilon_2
\left( x\right) =0 \ \ \ \ \ \forall x\in V$$
\end{prop}
\begin{proof}
Straightforward, (see \cite{Ka73}).
\end{proof}
\subsection{Construction of 2-Associative Bialgebras}
\

We construct an $(n+1)$-dimensional 2-associative bialgebra structure from an arbitrary  $n$-dimensional associative algebras.

 \begin{lemma} Let $\mathcal{A}=(V,\mu,\eta )$ be any unital
associative algebra. The  bialgebra
$\mathcal{K}_{1}\left( \mathcal{A}\right)
=(\widetilde{V},\mu_1,\eta_1,\Delta_1,\varepsilon_1 )$ is a unital infinitesimal
bialgebra.
 \end{lemma}
\begin{proof}
Since $\mathcal{K}_{1}\left( \mathcal{A}\right)$ is a bialgebra (prop.
\ref{kap1}), one has to check only the unital infinitesimal condition.
 Let $x,y\in V$ and $\tau$ be the usual flip. We have\newline
$\left( \mu_{1}\otimes \mu_{1}\right) \left( id_{V}\otimes
\tau\otimes id_{V}\right)  \left( i_{1}\otimes \Delta
_{1}+\Delta _{1}\otimes i_{2}-i_{1}\otimes i_{2}\right) \left(
x\otimes y\right) $

$=\left( \mu_{1}\otimes \mu_{1}\right) \left( id_{V}\otimes
\tau\otimes id_{V}\right) $

$\left( \left( x\otimes e_{1}\right) \otimes \left( y\otimes
e_{1}+e_{1}\otimes y-e_{2}\otimes y\right) +\left( x\otimes
e_{1}+e_{1}\otimes x-e_{2}\otimes x\right) \otimes \left(
e_{1}\otimes y \right) -x\otimes y\right) $

$=\mu\left( x\otimes y\right) \otimes e_{1}+e_{1}\otimes \mu%
\left( x\otimes y\right) -e_{2}\otimes \mu\left( x\otimes y\right) $

$=\Delta _{1}\left( \mu_{1}\left( x\otimes y\right) \right).
$\newline With similar computation, one gets  :\newline $\left(
\mu_{1}\otimes \mu_{1}\right) \left( id_{V}\otimes \tau\otimes
id_{V}\right)  \left( i_{1}\otimes \Delta _{1}+\Delta
_{1}\otimes i_{2}-i_{1}\otimes i_{2}\right) \left( e_{1}\otimes
x\right) $

$=\Delta _{1}\left( \mu_{1}\left( e_{1}\otimes x\right) \right)
=\Delta _{1}\left( x\right), $
\newline $\left( \mu_{1}\otimes
\mu_{1}\right) \left( id_{V}\otimes \tau\otimes id_{V}\right)
\left( i_{1}\otimes \Delta _{1}+\Delta _{1}\otimes
i_{2}-i_{1}\otimes i_{2}\right) \left( e_{1}\otimes e_{1}\right) $

$=\Delta _{1}\left( \mu_{1}\left( e_{1}\otimes e_{1}\right) \right)
=\Delta _{1}\left( e_{1}\right) .$\newline Then
$\mathcal{K}_{1}\left( A\right) $ is a unital infinitesimal bialgebra.
\end{proof}

\begin{rem} Let $\mathcal{A}=(V,\mu,\eta )$ be any unital
associative algebra. The bialgebra
$\mathcal{K}_{2}\left( \mathcal{A}\right) $ is not a unital infinitesimal
bialgebra, the unital infinitesimal condition is not fulfilled.
\end{rem}

\begin{rem}
Let $\A _2=(V,\mu_1,\mu_2,\eta)$ be a 2-associative algebra then we
have the same  coalgebra structure in the associated bialgebra
(resp. unital infinitesimal bialgebra) attached  to each unital associative algebra.
\end{rem}
  \begin{prop}
Let $\mathcal{A}=(V,\mu,\eta )$ and $\mathcal{A}'=(V,\mu ' ,\eta )$
be any two unital associative algebras over an $n$-dimensional
vector space $V$. Let $\mathcal{K}_{1}\left(
\mathcal{A}\right)=(\widetilde{V},\mu_1,\eta_1,\Delta_1,\varepsilon_1) $ and
$\mathcal{K}_{1}\left(
\mathcal{A}'\right)=(\widetilde{V},\mu_1',\eta_1,\Delta_1,\varepsilon_1) $ be
the associated bialgebras defined above. Then
$\mathfrak{B}_1=(\widetilde{V},\mu_1,\mu_1',\eta_1,\Delta_1,\varepsilon_1)$ is a
$(n+1)$-dimensional 2-associative bialgebra over the vector space
$\widetilde{V}=span(V,e_1)$ where $\eta_1(1)=e_1$.
\end{prop}
\begin{proof}
Since $\mathcal{K}_{1}\left(
\mathcal{A}\right)$ is a bialgebra and by the previous lemma
 $\mathcal{K}_{1}\left(
\mathcal{A}'\right)$ is an infinitesimal bialgebra, then
 $\mathfrak{B}_1=(\tilde{V},\mu_1,\mu_1',\eta _1,\Delta_1,\varepsilon_1)$ is 2-associative bialgebra.
\end{proof}

\begin{rem}
  Let $B=\left( V,\mu ,\eta,\Delta,\varepsilon \right) $ be a bialgebra, if the
comultiplication satisfies the unital infinitesimal condition \ref{inf} then
$\mathcal{B}2=\left( V,\mu ,\mu ,\eta,\Delta,\varepsilon \right) $ is a 2-associative
bialgebra.
\end{rem}

 \subsection{Construction of 2-Bialgebras}

\begin{prop}
 Let $V$ be an $n$-dimensional vector space over $\mathbb{K}$.
Let  $\mathcal{A}_{1}=(V,\mu _{1},\eta_1)$ and $\mathcal{A}_{2}=(V,\mu
_{2},\eta_2)$ be two unital associative algebras.

Let $\mathcal{K}_{i}\left( \A _{i}\right)=\left( \widetilde{V},\widetilde{\mu}_{i},\eta,\Delta
_{i},\varepsilon \right)\ \ i=1,2$ be the associated
bialgebras defined above.

Then
$$
\mathfrak{B}_{1}=\left( \widetilde{V},\widetilde{\mu}_{1},\widetilde{\mu}_{2},\eta,\Delta
_{1},\Delta _{2},\varepsilon \right)  \ \ \text{and} \ \
\mathfrak{B}_{2}=\left( \widetilde{V},\tilde{\mu}_{1},\widetilde{\mu}_{2},\eta,\Delta
_{1}^{cop},\Delta _{2},\varepsilon \right)
$$
are two $(n+1)$-dimensional 2-bialgebras on
$\widetilde{V}=span(V,e_1)$, where $e_1=\eta (1).$
\end{prop}
\begin{proof}
We show that $\mathfrak{B}_{1}$ is a 2-bialgebra. Since
$\mathfrak{B}_{i}\left( \A _{i}\right) $ $i=1,2$ are two  bialgebras, one has to prove the compatibility between $\widetilde{\mu}%
_{1}$ and $\Delta _{2}$, then  $\widetilde{\mu}_{2}$ and $\Delta
_{1}.$

Let $x,y\in V$.

 $\bullet \left(
\widetilde{\mu}_{1}\otimes \widetilde{\mu}_{1}\right) \circ \left(
id_{V}\otimes \tau\otimes id_{V}\right) \circ \left( \Delta
_{2}\otimes \Delta
_{2}\right) \left( x\otimes y\right) =\left( \widetilde{\mu}_{1}\otimes \tilde{%
\mu}_{1}\right) \circ \left( id_{V}\otimes \tau\otimes id_{V}\right)
$

$\left( \Delta _{2}\left( x\right) \otimes \Delta _{2}\left(
y\right) \right) $

$=\left( \widetilde{\mu}_{1}\otimes \widetilde{\mu}_{1}\right) \circ \left(
id_{V}\otimes \tau\otimes id_{V}\right) \left( \left(
e_{1}-e_{2}\right) \otimes x+x\otimes \left( e_{1}-e_{2}\right)
\right) \otimes $

$\left( \left( e_{1}-e_{2}\right) \otimes y+y\otimes \left(
e_{1}-e_{2}\right) \right) \ $

$=\left( e_{1}-e_{2}\right) \otimes \widetilde{\mu}_{1}\left( x\otimes y\right) +%
\widetilde{\mu}_{1}\left( x\otimes y\right) \otimes \left(
e_{1}-e_{2}\right) $

$=\Delta _{2}\circ \widetilde{\mu}_{1}\left( x\otimes y\right). $

Also

$\bullet \left( \widetilde{\mu}_{1}\otimes
\widetilde{\mu}_{1}\right) \circ \left( id_{V}\otimes \tau\otimes
id_{V}\right) \circ \left( \Delta _{2}\otimes \Delta
_{2}\right) \left( e_{2}\otimes e_{2}\right) =\Delta _{2}\circ \widetilde{\mu}%
_{1}\left( e_{2}\otimes e_{2}\right). $

$\bullet \left(
\widetilde{\mu}_{1}\otimes \widetilde{\mu}_{1}\right) \circ \left(
id_{V}\otimes \tau\otimes id_{V}\right) \circ \left( \Delta
_{2}\otimes \Delta
_{2}\right) \left( e_{2}\otimes x \right) =\Delta _{2}\circ \widetilde{\mu}_{1}\left( e_{2}\otimes x\right). $

$\bullet $ $\left(
\widetilde{\mu}_{1}\otimes \widetilde{\mu}_{1}\right) \circ \left(
id_{V}\otimes \tau\otimes id_{V}\right) \circ \left( \Delta
_{2}\otimes
\Delta _{2}\right) \left( x\otimes e_{2}\right) =\Delta _{2}\circ \widetilde{\mu}%
_{1}\left( x\otimes e_{2}\right) .$

 $\bullet
\left( \widetilde{\mu}_{2}\otimes \widetilde{\mu}_{2}\right) \circ \left(
id_{V}\otimes \tau\otimes id_{V}\right) \circ \left( \Delta
_{1}\otimes \Delta _{1}\right) \left( x\otimes y\right) $

$=\left( \widetilde{\mu}_{2}\otimes \widetilde{\mu}_{2}\right) \circ \left(
id_{V}\otimes \tau\otimes id_{V}\right) \left( \Delta _{1}\left(
x\right) \otimes \Delta _{1}\left( y\right) \right) $

$=\left( \widetilde{\mu}_{2}\otimes \widetilde{\mu}_{2}\right) \circ \left(
id_{V}\otimes \tau\otimes id_{V}\right) (\left( e_{1}\otimes
x+x\otimes e_{1}-e_{2}\otimes x\right) \otimes $

$\left( e_{1}\otimes y+y\otimes e_{1}-e_{2}\otimes y\right) )$

$=e_{1}\otimes \widetilde{\mu}_{2}\left( x\otimes y\right) +\widetilde{\mu}%
_{2}\left( x\otimes y\right) \otimes e_{1}-e_{2}\otimes \widetilde{\mu}%
_{2}\left( x\otimes y\right) $

$=\Delta _{1}\circ \widetilde{\mu}_{2}\left( x\otimes y\right) $.

One has also

 $\bullet \left( \widetilde{\mu}_{2}\otimes
\widetilde{\mu}_{2}\right) \circ \left( id_{V}\otimes \tau\otimes
id_{V}\right) \circ \left( \Delta _{1}\otimes \Delta
_{1}\right) \left( x\otimes e_{2}\right) =\Delta _{1}\circ \widetilde{\mu}%
_{2}\left( x\otimes e_{2}\right) $

$\bullet \left(
\widetilde{\mu}_{2}\otimes \widetilde{\mu}_{2}\right) \circ \left(
id_{V}\otimes \tau\otimes id_{V}\right) \circ \left( \Delta
_{1}\otimes \Delta
_{1}\right) \left( e_{2}\otimes x\right) =\Delta _{2}\circ \widetilde{\mu}%
_{1}\left( e_{2}\otimes x\right) .$

Henceforth, the compatibility between  $\widetilde{\mu}_{2}$ and
$\Delta _{1}$.
  Therefore $\mathfrak{B}_{1}$ is a 2-bialgebra.\newline Similar proof shows
that $\mathfrak{B}_{2}$ is a 2-bialgebra.
\end{proof}
We have also the obvious following corollary which gives 2-2-bialgebra starting with any two unital associative algebras.
\begin{cor}
 Let $V$ be an $n$-dimensional vector space over $\mathbb{K}$.
Let  $\mathcal{A}_{1}=(V,\mu _{1},\eta_1)$ and $\mathcal{A}_{2}=(V,\mu
_{2},\eta_2)$ be two unital associative algebras.

Let $\mathcal{K}_{1}\left( \A _{i}\right)\ \ i=1,2$ be the associated
bialgebras defined above.

Then
$$
\mathfrak{B}_{1}=\left( \widetilde{V},\widetilde{\mu}_{1},\widetilde{\mu}_{2},\eta,\Delta
_{1},\Delta _{1},\varepsilon \right)
$$
is a two $(n+1)$-dimensional 2-2-bialgebras on
$\widetilde{V}=span(V,e_1)$, where $e_1=\eta (1).$
\end{cor}
\begin{proof}Straightforward.
\end{proof}
\section{Classification  in low dimensions}
In this section, we show that for fixed dimension $n$, the 2-associative
bialgebras, 2-bialgebras and 2-2-bialgebras are endowed with a structure
of algebraic variety and a natural structure transport action which
 describes the set of isomorphic algebras. Solving such systems of polynomial
  equations leads to classifications of such structures. We aim at classifying  2-associative bialgebras and  2-bialgebras
 of dimension 2 and 3. First,
we establish the classification of bialgebras and infinitesimal
bialgebras. Then, we describe and enumerate 2 and 3-dimensional 2-associative bialgebras and 2-bialgebras.

Let $V$ be an $n$-dimensional vector space over $\mathbb{K}$.
Setting a basis $\left\{ e_{i}\right\} _{i=\left\{ 1,...,n\right\}
}$ of $V$,
 a multiplication $\mu $ (resp. a
comultiplication $\Delta $) is identified with its $n^{3}$ structure
constants $C_{ij}^{k}\in\mathbb{K}$ (resp. $D_{i}^{jk}$), where $\mu
\left( e_{i}\otimes e_{j}\right) =\sum_{k=1}^{n}C_{ij}^{k}e_{k}$ and
$\Delta \left( e_{i}\right) =\sum_{j,k=1}^{n}D_{i}^{jk}e_{j}\otimes
e_{k} $. The counit $\varepsilon $ is identified to its $n$
structure constants $\xi _{i}$. We assume that $e_{1}$ is the unit.

 A collection  $\{(
C_{ij}^{k},\tilde{C}_{ij}^{k},D_{i}^{jk},\xi _{i}),...i,j,k\in
\left\{ 1,...,n\right\} \}$ represents a 2-associative bialgebra if
the underlying multiplications, comultiplication, and the counit
satisfy the appropriate conditions which translate to following
polynomial equations.

$\left( A_{1}\right) \left\{
\begin{array}{c}
\sum_{l=1}^{n}\left(
C_{ij}^{l}C_{lk}^{s}-C_{jk}^{l}C_{il}^{s}\right) =0\ \
\  \\
C_{1i}^{j}=C_{i1}^{j}=\delta _{ij}\ \ \ \ \ \ \ \ \ \ \ \ \ \ \ \ \
\end{array}
\right. \ \ \ \ \ \ \ \ \ \ \ \ \ \ \ \ \ \ \ \ \ \ \forall
i,j,k,s\in \left\{ 1,...,n\right\} $

$\left( A_{2}\right) \left\{
\begin{array}{c}
\sum_{l=1}^{n}\left( \tilde{C}_{ij}^{l}\tilde{C}_{lk}^{s}-\tilde{C}_{jk}^{l}%
\tilde{C}_{il}^{s}\right) =0\  \\
\tilde{C}_{1i}^{j}=\tilde{C}_{i1}^{j}=\delta _{ij}\ \ \ \ \ \ \ \ \
\ \ \ \ \ \ \
\end{array}
\right. \ \ \ \ \ \ \ \ \ \ \ \ \ \ \ \ \ \ \ \forall i,j,k,s\in
\left\{ 1,...,n\right\} $

$\left( A_{3}\right) \left\{
\begin{array}{c}
\sum_{l=1}^{n}\left(
D_{s}^{lk}D_{l}^{ij}-D_{s}^{il}D_{l}^{jk}\right) =0\
\\
\sum_{l=1}^{n}D_{i}^{jl}\zeta _{l}=\sum_{l=1}^{n}D_{i}^{lj}\zeta
_{l}=\delta _{ij}\ \ \ \ \
\end{array}
\right. \ \ \ \ \ \ \ \ \ \ \ \ \forall i,j,k,s\in \left\{
1,...,n\right\} $

$\left( A_{4}\right) \left\{
\begin{array}{c}
\sum_{l=1}^{n}\tilde{C}_{ij}^{l}D_{l}^{ks}-%
\sum_{r,t,p,q=1}^{n}D_{i}^{rt}D_{j}^{pq}\tilde{C}_{rp}^{k}\tilde{C}%
_{tq}^{s}=0 \\
D_{1}^{11}=1,D_{1}^{ij}=0\text{ \ \ }\left( i,j\right) \neq \left(
1,1\right)  \\
\zeta _{1}=1,\sum_{l=1}^{n}\tilde{C}_{ij}^{l}\zeta _{l}=\zeta
_{i}\zeta _{j}
\end{array}
\right. \ \forall i,j,k,s\in \left\{ 1,...,n\right\} $

$\left( A_{5}\right) \left\{
\begin{array}{c}
\sum_{l=1}^{n}\left(
D_{j}^{lj}C_{il}^{i}+D_{i}^{il}C_{lj}^{j}-D_{l}^{ij}C_{ij}^{l}\right) =1%
\text{ \ \ \ \ \ \ \ \ \ \ \ \ \ \ \ \ \ \ \ \ \ \ \ \ \ \ } \\
\sum_{l=1}^{n}\left(
D_{j}^{lk}C_{il}^{s}+D_{i}^{sl}C_{lj}^{k}-D_{l}^{sk}C_{ij}^{l}\right) =0%
\text{\ }\left( i,j\right) \neq \left( s,k\right)
\end{array}
\right. \forall i,j,k,s\in \left\{ 1,...,n\right\} $

 Then, the set of $n$-dimensional 2-associative bialgebras,
which we denote by $2Ass\mathcal{B}_n$, carries a structure of algebraic
variety imbedded in $\mathbb{K}^{3n^3+n}$ with its natural structure
of algebraic variety.

  Similarly,  a collection  $\{ (
C_{ij}^{k},\tilde{C}_{ij}^{k},D_{i}^{jk},\tilde{D}_{i}^{jk},\xi
_{i},\tilde{\xi} _{i}),...i,j,k\in \left\{ 1,...,n\right\} \}$
represents a 2-bialgebra if it satisfies the following system

 $\left( B_{1}\right)
\left\{
\begin{array}{c}
\sum_{l=1}^{n}\left(
C_{ij}^{l}C_{lk}^{s}-C_{jk}^{l}C_{il}^{s}\right) =0\ \
\\
C_{1i}^{j}=C_{i1}^{j}=\delta _{ij}\ \ \ \ \ \ \ \ \ \ \ \ \ \ \
\end{array}
\right. $ \ \ \ \ \ \ \ \ $\ \ \forall i,j,k,s\in \left\{
1,...,n\right\} $

$\left( B_{2}\right) \left\{
\begin{array}{c}
\sum_{l=1}^{n}\left( \tilde{C}_{ij}^{l}\tilde{C}_{lk}^{s}-\tilde{C}_{jk}^{l}%
\tilde{C}_{il}^{s}\right) =0\ \  \\
\tilde{C}_{1i}^{j}=\tilde{C}_{i1}^{j}=\delta _{ij}\ \ \ \ \ \ \ \ \
\ \ \ \ \ \ \
\end{array}
\right. \ $\ $\ \ \ \ \ \ \ \ \ \ \ \forall i,j,k,s\in \left\{
1,...,n\right\} $

$\left( B_{3}\right) \left\{
\begin{array}{c}
\sum_{l=1}^{n}\left(
D_{s}^{lk}D_{l}^{ij}-D_{s}^{il}D_{l}^{jk}\right) =0\ \
\ \  \\
\sum_{l=1}^{n}D_{i}^{jl}\zeta _{l}=\stackrel{n}{\sum{l=1}{\sum }}%
D_{i}^{lj}\zeta _{l}=\delta _{ij}\ \ \
\end{array}
\right. $ $\ \ \ \ \ \ \ \ \forall i,j,k,s\in \left\{
1,...,n\right\} $

$\left( B_{4}\right) \left\{
\begin{array}{c}
\sum_{l=1}^{n}\left( \tilde{D}_{s}^{lk}\tilde{D}_{l}^{ij}-\tilde{D}_{s}^{il}%
\tilde{D}_{l}^{jk}\right) =0\ \ \  \\
\sum_{l=1}^{n}\tilde{D}_{i}^{jl}\tilde{\zeta}_{l}=\sum_{l=1}^{n}\tilde{D}%
_{i}^{lj}\tilde{\zeta}_{l}=\delta _{ij}\ \ \
\end{array}
\right. $ $\ \ \ \ \ \ \ \ \ \ \ \forall i,j,k,s\in \left\{
1,...,n\right\} $

$\left( B_{5}\right) \left\{
\begin{array}{c}
\sum_{l=1}^{n}C_{ij}^{l}D_{l}^{ks}-%
\sum_{r,t,p,q=1}^{n}D_{i}^{rt}D_{j}^{pq}C_{rp}^{k}C_{tq}^{s}=0 \\
D_{1}^{11}=1,D_{1}^{ij}=0\text{ \ \ \ \ \ \ \ \ \ \ }\left(
i,j\right) \neq
\left( 1,1\right)  \\
\zeta _{1}=1,\sum_{l=1}^{n}C_{ij}^{l}\zeta _{l}=\zeta _{i}\zeta
_{j}\text{ \ \ \ \ \ \ \ \ \ \ \ \ \ \ \ \ \ \ \ \ \ }
\end{array}
\right. \ \ \forall i,j,k,s\in \left\{ 1,...,n\right\} $

$\left( B_{6}\right) \left\{
\begin{array}{c}
\sum_{l=1}^{n}\tilde{C}_{ij}^{l}\tilde{D}_{l}^{ks}-\sum_{r,t,p,q=1}^{n}%
\tilde{D}_{i}^{rt}\tilde{D}_{j}^{pq}\tilde{C}_{rp}^{k}\tilde{C}_{tq}^{s}=0
\\
\tilde{D}_{1}^{11}=1,\tilde{D}_{1}^{ij}=0\text{ \ \ \ \ \ \ \ \ \ \ \ \ \ }%
\left( i,j\right) \neq \left( 1,1\right)  \\
\tilde{\zeta}_{1}=1,\sum_{l=1}^{n}\tilde{C}_{ij}^{l}\tilde{\zeta}_{l}=\tilde{%
\zeta}_{i}\tilde{\zeta}_{j}\text{ \ \ \ \ \ \ \ \ \ \ \ \ \ \ \ \ \
\ \ \ \ \ }
\end{array}
\right. \ \ \ \ \forall i,j,k,s\in \left\{ 1,...,n\right\} $

$\left( B_{7}\right) \left\{
\begin{array}{c}
\sum_{l=1}^{n}\tilde{C}_{ij}^{l}D_{l}^{ks}-%
\sum_{r,t,p,q=1}^{n}D_{i}^{rt}D_{j}^{pq}\tilde{C}_{rp}^{k}\tilde{C}%
_{tq}^{s}=0 \\
\sum_{l=1}^{n}C_{ij}^{l}\tilde{D}_{l}^{ks}-\sum_{r,t,p,q=1}^{n}\tilde{D}%
_{i}^{rt}\tilde{D}_{j}^{pq}C_{rp}^{k}C_{tq}^{s}=0
\end{array}
\right. \ \ \ \ \ \ \ \forall i,j,k,s\in \left\{ 1,...,n\right\} $

$\left( B_{8}\right) \left\{
\begin{array}{c}
\sum_{l=1}^{n}\tilde{C}_{ij}^{l}\zeta _{l}=\zeta _{i}\zeta _{j} \\
\sum_{l=1}^{n}C_{ij}^{l}\tilde{\zeta}_{l}=\tilde{\zeta}_{i}\tilde{\zeta}_{j}
\end{array}
\right. \ \ \ \ \ \ \ \ \ \ \ \ \ \ \ \ \ \ \ \ \ \ \ \ \ \ \ \ \ \
\ \ \ \ \forall i,j\in \left\{ 1,...,n\right\} $ $\ \ \ \ \ \ \ \ \
$

Then, the set of $n$-dimensional 2-bialgebras, which we denote by
$2\mathcal{B}_n$, carries a structure of algebraic variety imbedded in
$\mathbb{K}^{4n^3+2n}$.

Similarly, we have an algebraic variety structure on the set of $n$-dimensional 2-2-bialgebras, which we denote by  $22\mathcal{B}_n.$

The "structure transport" action is defined by the action of $GL_n(V)$
on $2Ass\mathcal{B}_n$ (similarly on  $2\mathcal{B}_n$ and $22\mathcal{B}_n$). It corresponds to the change of basis.

Let $\mathcal{B}=\left( V,\mu_1,\eta_1,\mu_2,\eta_2
 ,\Delta ,\varepsilon \right) $  be
a 2-associative bialgebras and   $f: V\rightarrow
V$ be an invertible endomorphism, then the action of $f$ on $\mathcal{B}$ transports the 2-associative bialgebra structure into a  2-associative bialgebra $\mathcal{B}'=\left( V,\mu_1',\eta_1',\mu_2',\eta_2',\Delta' ,\varepsilon' \right) $ defined by
\begin{align*}
\mu_1^{\prime}  =f\circ \mu_1\circ \left( f^{-1}\otimes f^{-1}\right) \quad \text {and} \quad
  \eta_1^{\prime }=f\circ \eta_1 \\
\mu_2^{\prime } =f\circ
\mu_2\circ \left( f^{-1}\otimes f^{-1}\right) \quad \text {and} \quad \eta_2^{\prime
}=f\circ \eta_2\\
 \Delta^{\prime }=(f\otimes f)\circ \Delta\circ f^{-1} \quad \text {and} \quad \varepsilon^{\prime } =\varepsilon \circ f^{-1}\mathit{\ }
\end{align*}

 \subsection{Classifications in dimension $2$ }

The set of 2-dimensional unital associative algebras yields two
non-isomorphic algebras (see \cite{Ga74}). Let $\{ e_1,e_2\}$ be a
basis of $\K ^2$, then the algebras are given by the following
non-trivial products.

$\bullet \mu _{1}^{2}\left( e_{1}\otimes e_{i}\right) =\mu _{1}^{2}\left(
e_{i}\otimes e_{1}\right) =e_{i},\ i=1,2,\ \mu _{1}^{2}\left(
e_{2}\otimes e_{2}\right) =e_{2}\medskip $

$\bullet \mu _{2}^{2}\left( e_{1}\otimes e_{i}\right) =\mu _{2}^{2}\left(
e_{i}\otimes e_{1}\right) =e_{i},\ i=1,2,\ \mu _{1}^{2}\left(
e_{2}\otimes e_{2}\right) =0\medskip $

In the sequel we consider that all the algebras are unital and the unit
$\eta$  corresponds to $e_1$.

In the following, we list the coalgebras which, combined with
$\mu_1$, give  bialgebra structures  (up to isomorphism).

 $\bullet \Delta _{1,1}^{2}\left( e_{1}\right) =e_{1}\otimes e_{1},\ \Delta
_{1,1}^{2}\left( e_{2}\right) =e_{1}\otimes e_{2}+e_{2}\otimes
e_{1}-2e_{2}\otimes e_{2},$

$\varepsilon _{1,1}^{2}\left( e_{1}\right) =1,\ \varepsilon
_{1,1}^{2}\left( e_{2}\right) =0.$

 $\bullet \Delta _{1,2}^{2}\left( e_{1}\right) =e_{1}\otimes
e_{1},\ \Delta _{1,2}^{2}\left( e_{2}\right) =e_{2}\otimes e_{2},$

$ \varepsilon _{1,2}^{2}\left( e_{1}\right) =1,\ \varepsilon
_{1,2}^{2}\left( e_{2}\right) =1.$

  $\bullet \Delta _{1,3}^{2}\left( e_{1}\right) =e_{1}\otimes
e_{1};\Delta _{1,3}^{2}\left( e_{2}\right) =e_{1}\otimes
e_{2}+e_{2}\otimes e_{1}-e_{2}\otimes e_{2},$

  $\varepsilon
_{1,3}^{2}\left( e_{1}\right) =1,\ \varepsilon _{1,3}^{2}\left(
e_{2}\right) =0.$

We cannot associate to $\mu _{2}^{2}$ neither a bialgebra structure
nor an infinitesimal bialgebra structure.

Direct calculations show that there is only one bialgebra which carries also a structure of infinitesimal bialgebra. This leads to the following classification of 2-dimensional 2-associative bialgebras.
\begin{prop}
Every 2-dimensional 2-associative bialgebra is isomorphic to $\left(
\mathbb{K}^2,\mu _{1}^{2},\mu _{1}^{2},\eta, \Delta _{1,2}^{2},\varepsilon
_{1,2}^{2}\right) .$
\end{prop}

 The 2-dimensional 2-bialgebras are given by the following
 proposition.
 \begin{prop}

Every 2-bialgebra of type (1,1) is isomorphic to one of the
following 2-bialgebras
$$\left( V,\mu _{1}^{2},\mu _{1}^{2},\eta,\Delta _{1,1}^{2},\Delta
_{1,1}^{2},\varepsilon _{1,1}^{2}\right) ,\ \left( V,\mu _{1}^{2},\mu
_{1}^{2},\eta,\Delta _{1,2}^{2},\Delta _{1,2}^{2},\varepsilon
_{1,2}^{2}\right),
 \  \left( V,\mu _{1}^{2},\mu _{1}^{2},\eta,\Delta _{1,3}^{2},\Delta
_{1,3}^{2},\varepsilon _{1,3}^{2}\right) .$$
Every 2-bialgebra of type (1,2) is isomorphic to one of the
following 2-bialgebras
$$\left( V,\mu _{1}^{2},\mu _{1}^{2},\eta,\Delta _{1,1}^{2},\Delta
_{1,2}^{2},\varepsilon _{1,2}^{2}\right) ,\ \left( V,\mu _{1}^{2},\mu
_{1}^{2},\eta,\Delta _{1,1}^{2},\Delta _{1,3}^{2},\varepsilon
_{1,1}^{2}\right), \ \left( V,\mu _{1}^{2},\mu _{1}^{2},\eta,\Delta _{1,2}^{2},\Delta
_{1,3}^{2},\varepsilon _{1,3}^{2}\right) .$$
There is no 2-bialgebras of type (2,2) and (2,1) in dimension 2.
\end{prop}
\begin{rem}
There is only one 2-dimensional 2-2-bialgebra which is given by $$\left(
\mathbb{K}^2,\mu _{1}^{2},\mu _{1}^{2},\eta, \Delta _{1,2}^{2},\Delta _{1,2}^{2},\varepsilon
_{1,2}^{2}\right) .$$
\end{rem}

 \subsection{Classifications in dimension $3$ }

First, we recall the classification of 3-dimensional unital
associative algebras (see \cite{Ga74}). Let $\{ e_1,e_2,e_3\}$ be a basis of $\K ^3$, then the algebras are given by the following non-trivial products.

$\bullet \mu _{1}^{3}\left( e_{1}\otimes e_{i}\right) =\mu _{1}^{3}\left(
e_{i}\otimes e_{1}\right) =e_{i} \ \ i=1,2,3,$

$\mu _{1}^{3}\left( e_{j}\otimes e_{2}\right) =\mu _{1}^{3}\left(
e_{2}\otimes e_{j}\right) =e_{j}\ \ j=2,3,\ \mu _{1}^{3}\left(
e_{3}\otimes e_{3}\right) =e_{3}.$

$\bullet \mu _{2}^{3}\left( e_{1}\otimes e_{i}\right) =\mu _{2}^{3}\left(
e_{i}\otimes e_{1}\right) =e_{i} \ \ i=1,2,3,$

$\mu _{2}^{3}\left( e_{j}\otimes e_{2}\right) =\mu _{2}^{3}\left(
e_{2}\otimes e_{j}\right) =e_{j}\ \ j=2,3,~\mu _{2}^{3}\left(
e_{3}\otimes e_{3}\right) =0.$

 $\bullet \mu _{3}^{3}\left( e_{1}\otimes e_{i}\right) =\mu _{3}^{3}\left(
e_{i}\otimes e_{1}\right) =e_{i} \ \ i=1,2,3,$ $ \mu _{3}^{3}\left(
e_{2}\otimes e_{2}\right) =e_{2}.$

$\bullet \mu _{4}^{3}\left( e_{1}\otimes e_{i}\right) =\mu _{4}^{3}\left(
e_{i}\otimes e_{1}\right) =e_{i},\quad i=1,2,3.$

$\bullet \mu _{5}^{3}\left( e_{1}\otimes e_{i}\right) =\mu _{5}^{3}\left(
e_{i}\otimes e_{1}\right) =e_{i}\ \ i=1,2,3;\ \ \mu _{5}^{3}\left(
e_{2}\otimes e_{j}\right) =e_{j} \ \ j=2,3.$

\vspace{0.3cm}
Thanks to computer algebra, we obtain the following coalgebras
associated to the previous algebras in order to obtain a bialgebra
structures.
We denote the comultiplications by $\Delta _{i,j}^{3}$ and the counits by $%
\varepsilon _{i,j}^{3}$, where $i$ indicates the item of the
multiplication and $j$ the item of the comultiplication which combined with the multiplication $i$ determine a bialgebra.

\vspace{0.3cm}
For the multiplication $\mu _{1}^{3}$, we have

\begin{enumerate}
\item  $ \Delta _{1,1}^{3}\left( e_{1}\right) =e_{1}\otimes
e_{1};\ \Delta _{1,1}^{3}\left( e_{2}\right) =e_{1}\otimes
e_{2}+e_{2}\otimes e_{1}-e_{2}\otimes e_{2};\ \Delta _{1,1}^{3}\left(
e_{3}\right) =e_{1}\otimes e_{3}+e_{3}\otimes e_{1}-2e_{3}\otimes
e_{3};\ \varepsilon _{1,1}^{3}\left( e_{1}\right) =1;\ \varepsilon
_{1,1}^{3}\left( e_{2}\right) =0;\ \varepsilon _{1,1}^{3}\left(
e_{3}\right) =0.$

\item  $ \Delta _{1,2}^{3}\left( e_{1}\right) =e_{1}\otimes
e_{1};\ \Delta _{1,2}^{3}\left( e_{2}\right) =e_{1}\otimes
e_{2}+e_{2}\otimes e_{1}-e_{2}\otimes e_{2};\ \Delta _{1,2}^{3}\left(
e_{3}\right) =e_{1}\otimes e_{3}+e_{3}\otimes e_{1}-e_{3}\otimes
e_{3};\ \varepsilon _{1,2}^{3}\left( e_{1}\right) =1;\ \varepsilon
_{1,2}^{3}\left( e_{2}\right) =0;\ \varepsilon _{1,2}^{3}\left(
e_{3}\right) =0.$

\item  $ \Delta _{1,3}^{3}\left( e_{1}\right) =e_{1}\otimes
e_{1};\ \Delta _{1,3}^{3}\left( e_{2}\right) =e_{1}\otimes
e_{2}+e_{2}\otimes e_{1}-e_{2}\otimes e_{2};\ \Delta _{1,3}^{3}\left(
e_{3}\right) =e_{1}\otimes e_{3}-e_{2}\otimes e_{3}+e_{3}\otimes
e_{1}-e_{3}\otimes e_{2}-e_{3}\otimes e_{3};\ \varepsilon
_{1,3}^{3}\left( e_{1}\right) =1;\ \varepsilon _{1,3}^{3}\left(
e_{2}\right) =0;\ \varepsilon _{1,3}^{3}\left( e_{3}\right) =0. $

\item  $ \Delta _{1,4}^{3}\left( e_{1}\right) =e_{1}\otimes
e_{1};\ \Delta _{1,4}^{3}\left( e_{2}\right) =e_{1}\otimes
e_{2}+e_{2}\otimes e_{1}-e_{2}\otimes e_{2};\ \Delta _{1,4}^{3}\left(
e_{3}\right) =e_{1}\otimes e_{3}-e_{2}\otimes e_{3}+e_{3}\otimes
e_{1}-e_{3}\otimes e_{2};\ \varepsilon _{1,4}^{3}\left( e_{1}\right)
=1;\ \varepsilon _{1,4}^{3}\left( e_{2}\right) =0;\ \varepsilon
_{1,4}^{3}\left( e_{3}\right) =0$

\item  $ \Delta _{1,5}^{3}\left( e_{1}\right) =e_{1}\otimes
e_{1};\ \Delta _{1,5}^{3}\left( e_{2}\right) =e_{1}\otimes
e_{2}+e_{2}\otimes e_{1}-e_{2}\otimes e_{2};\ \Delta _{1,5}^{3}\left(
e_{3}\right) =e_{1}\otimes e_{3}+e_{3}\otimes e_{1}-e_{2}\otimes
e_{3};\ \varepsilon _{1,5}^{3}\left( e_{1}\right) =1;\ \varepsilon
_{1,5}^{3}\left( e_{2}\right) =0;\ \varepsilon _{1,5}^{3}\left(
e_{3}\right) =0$

\item  $ \Delta _{1,6}^{3}\left( e_{1}\right) =e_{1}\otimes
e_{1};\ \Delta _{1,6}^{3}\left( e_{2}\right) =e_{1}\otimes
e_{2}+e_{2}\otimes e_{1}-e_{2}\otimes e_{2};\ \Delta _{1,6}^{3}\left(
e_{3}\right) =e_{1}\otimes e_{3}+e_{3}\otimes e_{1}-e_{3}\otimes
e_{2};\ \varepsilon _{1,6}^{3}\left( e_{1}\right) =1;\ \varepsilon
_{1,6}^{3}\left( e_{2}\right) =0;\ \varepsilon _{1,6}^{3}\left(
e_{3}\right) =0.$

\item  $ \Delta _{1,7}^{3}\left( e_{1}\right) =e_{1}\otimes
e_{1};\ \Delta _{1,7}^{3}\left( e_{2}\right) =e_{2}\otimes
e_{2};\ \Delta _{1,7}^{3}\left( e_{3}\right) =e_{2}\otimes
e_{3}+e_{3}\otimes e_{2}-2e_{3}\otimes e_{3};\ \varepsilon
_{1,7}^{3}\left( e_{1}\right) =1;\ \varepsilon _{1,7}^{3}\left(
e_{2}\right) =1;\ \varepsilon _{1,7}^{3}\left( e_{3}\right) =0.$

\item  $ \Delta _{1,8}^{3}\left( e_{1}\right) =e_{1}\otimes
e_{1};\ \Delta _{1,8}^{3}\left( e_{2}\right) =e_{2}\otimes
e_{2};\ \Delta _{1,8}^{3}\left( e_{3}\right) =e_{2}\otimes
e_{3}+e_{3}\otimes e_{2}-e_{3}\otimes e_{3};\ \varepsilon
_{1,8}^{3}\left( e_{1}\right) =1;\ \varepsilon _{1,8}^{3}\left(
e_{2}\right) =1;\ \varepsilon _{1,8}^{3}\left( e_{3}\right) =0.$

\item  $ \Delta _{1,9}^{3}\left( e_{1}\right) =e_{1}\otimes
e_{1};\ \Delta _{1,9}^{3}\left( e_{2}\right) =e_{1}\otimes
e_{3}+e_{2}\otimes e_{2}-e_{2}\otimes e_{3}+e_{3}\otimes
e_{1}-e_{3}\otimes e_{2};\ \Delta _{1,9}^{3}\left( e_{3}\right)
=e_{1}\otimes e_{3}+e_{3}\otimes e_{1}-e_{3}\otimes
e_{3};\ \varepsilon _{1,9}^{3}\left( e_{1}\right) =1;\ \varepsilon
_{1,9}^{3}\left( e_{2}\right) =1;\ \varepsilon _{1,9}^{3}\left(
e_{3}\right) =0.$

\item  $ \Delta _{1,10}^{3}\left( e_{1}\right) =e_{1}\otimes
e_{1};\ \Delta _{1,10}^{3}\left( e_{2}\right) =e_{1}\otimes
e_{3}+e_{2}\otimes e_{2}-e_{2}\otimes e_{3}+e_{3}\otimes
e_{1}-e_{3}\otimes e_{2}+e_{3}\otimes e_{3};\ \Delta _{1,10}^{3}\left(
e_{3}\right) =e_{1}\otimes e_{3}+e_{3}\otimes e_{1}-2e_{3}\otimes
e_{3};\ \varepsilon _{1,10}^{3}\left( e_{1}\right) =1;\ \varepsilon
_{1,10}^{3}\left( e_{2}\right) =1;\ \varepsilon _{1,10}^{3}\left(
e_{3}\right) =0.$

\item  $ \Delta _{1,11}^{3}\left( e_{1}\right) =e_{1}\otimes
e_{1};\ \Delta _{1,11}^{3}\left( e_{2}\right) =e_{2}\otimes
e_{2}+e_{3}\otimes e_{1}-e_{3}\otimes e_{2};\ \Delta _{1,11}^{3}\left(
e_{3}\right) =e_{2}\otimes e_{3}+e_{3}\otimes e_{1}-e_{3}\otimes
e_{3};\ \varepsilon _{1,11}^{3}\left( e_{1}\right) =1;\ \varepsilon
_{1,11}^{3}\left( e_{2}\right) =1;\ \varepsilon _{1,11}^{3}\left(
e_{3}\right) =0.$

\item  $ \Delta _{1,12}^{3}\left( e_{1}\right) =e_{1}\otimes
e_{1};\ \Delta _{1,12}^{3}\left( e_{2}\right) =e_{1}\otimes
e_{3}+e_{2}\otimes e_{2}-e_{2}\otimes e_{3};\ \Delta _{1,12}^{3}\left(
e_{3}\right) =e_{1}\otimes e_{3}+e_{3}\otimes e_{2}-e_{3}\otimes
e_{3};\ \varepsilon _{1,12}^{3}\left( e_{1}\right) =1;\ \varepsilon
_{1,12}^{3}\left( e_{2}\right) =1;\ \varepsilon _{1,12}^{3}\left(
e_{3}\right) =0.$

\item  $ \Delta _{1,13}^{3}\left( e_{1}\right) =e_{1}\otimes
e_{1};\ \Delta _{1,13}^{3}\left( e_{2}\right) =e_{1}\otimes
e_{2}-e_{1}\otimes e_{3}+e_{2}\otimes e_{1}-2e_{2}\otimes
e_{2}+2e_{2}\otimes e_{3}-e_{3}\otimes e_{1}+2e_{3}\otimes
e_{2}-e_{3}\otimes e_{3};\ \Delta _{1,13}^{3}\left( e_{3}\right)
=e_{2}\otimes e_{3}+e_{3}\otimes e_{2}-2e_{3}\otimes
e_{3};\ \varepsilon _{1,13}^{3}\left( e_{1}\right) =1;\ \varepsilon
_{1,13}^{3}\left( e_{2}\right) =1;\ \varepsilon _{1,13}^{3}\left(
e_{3}\right) =1.$

\item  $ \Delta _{1,14}^{3}\left( e_{1}\right) =e_{1}\otimes
e_{1};\ \Delta _{1,14}^{3}\left( e_{2}\right) =e_{1}\otimes
e_{2}-e_{1}\otimes e_{3}+e_{2}\otimes e_{1}-e_{2}\otimes
e_{2}+e_{2}\otimes e_{3}-e_{3}\otimes e_{1}+e_{3}\otimes
e_{2};\ \Delta _{1,14}^{3}\left( e_{3}\right) =e_{2}\otimes
e_{3}+e_{3}\otimes e_{2}-e_{3}\otimes e_{3};\ \varepsilon
_{1,14}^{3}\left( e_{1}\right) =1;\ \varepsilon _{1,14}^{3}\left(
e_{2}\right) =1;\ \varepsilon _{1,14}^{3}\left( e_{3}\right) =1.$

\item  $ \Delta _{1,15}^{3}\left( e_{1}\right) =e_{1}\otimes
e_{1};\ \Delta _{1,15}^{3}\left( e_{2}\right) =e_{2}\otimes
e_{2};\ \Delta _{1,15}^{3}\left( e_{3}\right) =e_{3}\otimes
e_{3};\ \varepsilon _{1,15}^{3}\left( e_{1}\right) =1;\ \varepsilon
_{1,15}^{3}\left( e_{2}\right) =1;\ \varepsilon _{1,15}^{3}\left(
e_{3}\right) =1.$

\item  $ \Delta _{1,16}^{3}\left( e_{1}\right) =e_{1}\otimes
e_{1};\ \Delta _{1,16}^{3}\left( e_{2}\right) =e_{2}\otimes
e_{2};\ \Delta _{1,16}^{3}\left( e_{3}\right) =e_{2}\otimes
e_{2}-e_{2}\otimes e_{3}-e_{3}\otimes e_{2}+2e_{3}\otimes
e_{3};\ \varepsilon _{1,16}^{3}\left( e_{1}\right) =1;\ \varepsilon
_{1,16}^{3}\left( e_{2}\right) =1;\ \varepsilon _{1,16}^{3}\left(
e_{3}\right) =1.$

\item  $ \Delta _{1,17}^{3}\left( e_{1}\right) =e_{1}\otimes
e_{1};\ \Delta _{1,17}^{3}\left( e_{2}\right) =e_{2}\otimes
e_{3}+e_{3}\otimes e_{2}-e_{3}\otimes e_{3};\ \Delta _{1,17}^{3}\left(
e_{3}\right) =e_{3}\otimes e_{3};\ \varepsilon _{1,17}^{3}\left(
e_{1}\right) =1;\ \varepsilon _{1,17}^{3}\left( e_{2}\right)
=1;\ \varepsilon _{1,17}^{3}\left( e_{3}\right) =1.$

\item  $ \Delta _{1,18}^{3}\left( e_{1}\right) =e_{1}\otimes
e_{1};\ \Delta _{1,18}^{3}\left( e_{2}\right) =e_{2}\otimes
e_{1}-e_{3}\otimes e_{1}+e_{3}\otimes e_{2};\ \Delta _{1,18}^{3}\left(
e_{3}\right) =e_{3}\otimes e_{3};\ \varepsilon _{1,18}^{3}\left(
e_{1}\right) =1;\ \varepsilon _{1,18}^{3}\left( e_{2}\right)
=1;\ \varepsilon _{1,18}^{3}\left( e_{3}\right) =1.$
\end{enumerate}

For the multiplication $\mu _{2}^{3}$, we have

\begin{enumerate}
\item $ \Delta _{2,1}^{3}\left( e_{1}\right) =e_{1}\otimes
e_{1};\ \Delta _{2,1}^{3}\left( e_{2}\right) =e_{1}\otimes
e_{2}+e_{2}\otimes e_{1}-e_{2}\otimes e_{2};\ \Delta _{2,1}^{3}\left(
e_{3}\right) =e_{1}\otimes e_{3}+e_{3}\otimes e_{1}-e_{3}\otimes
e_{2};\ \varepsilon _{2,1}^{3}\left( e_{1}\right) =1;\ \varepsilon
_{2,1}^{3}\left( e_{2}\right) =0;\ \varepsilon _{2,1}^{3}\left(
e_{3}\right) =0.$

\item  $ \Delta _{2,2}^{3}\left( e_{1}\right) =e_{1}\otimes
e_{1};\ \Delta _{2,2}^{3}\left( e_{2}\right) =e_{1}\otimes
e_{2}+e_{2}\otimes e_{1}-e_{2}\otimes e_{2};\ \Delta _{2,2}^{3}\left(
e_{3}\right) =e_{1}\otimes e_{3}+e_{2}\otimes e_{3}+e_{3}\otimes
e_{1};\ \varepsilon _{2,2}^{3}\left( e_{1}\right) =1;\ \varepsilon
_{2,2}^{3}\left( e_{2}\right) =0;\ \varepsilon _{2,2}^{3}\left(
e_{3}\right) =0.$

\item  $ \Delta _{2,3}^{3}\left( e_{1}\right) =e_{1}\otimes
e_{1};\ \Delta _{2,3}^{3}\left( e_{2}\right) =e_{1}\otimes
e_{2}+e_{2}\otimes e_{1}-e_{2}\otimes e_{2};\ \Delta _{2,3}^{3}\left(
e_{3}\right) =e_{1}\otimes e_{3}-e_{2}\otimes e_{3}+e_{3}\otimes
e_{1}-e_{3}\otimes e_{2}+\lambda e_{3}\otimes e_{3};\ \varepsilon
_{2,3}^{3}\left( e_{1}\right) =1;\ \varepsilon _{2,3}^{3}\left(
e_{2}\right) =0;\ \varepsilon _{2,3}^{3}\left( e_{3}\right) =0.$
\end{enumerate}

For the multiplication $\mu _{3}^{3}$, we have

\begin{enumerate}
\item  $ \Delta _{3,1}^{3}\left( e_{1}\right) =e_{1}\otimes
e_{1};\ \Delta _{3,1}^{3}\left( e_{2}\right) =e_{2}\otimes
e_{2};\ \Delta _{3,1}^{3}\left( e_{3}\right) =e_{2}\otimes
e_{3}+e_{3}\otimes e_{2};\ \varepsilon _{3,1}^{3}\left( e_{1}\right)
=1;\ \varepsilon _{3,1}^{3}\left( e_{2}\right) =1;\ \varepsilon
_{3,1}^{3}\left( e_{3}\right) =0 .$

\item  $ \Delta _{3,2}^{3}\left( e_{1}\right) =e_{1}\otimes
e_{1};\ \Delta _{3,2}^{3}\left( e_{2}\right) =e_{2}\otimes
e_{2};\ \Delta _{3,2}^{3}\left( e_{3}\right) =e_{1}\otimes
e_{3}+e_{3}\otimes e_{2};\ \varepsilon _{3,2}^{3}\left( e_{1}\right)
=1;\ \varepsilon _{3,2}^{3}\left( e_{2}\right) =1;\ \varepsilon
_{3,2}^{3}\left( e_{3}\right) =0. $

\item  $ \Delta _{3,3}^{3}\left( e_{1}\right) =e_{1}\otimes
e_{1};\ \Delta _{3,3}^{3}\left( e_{2}\right) =e_{2}\otimes
e_{2};\ \Delta _{3,3}^{3}\left( e_{3}\right) =e_{2}\otimes
e_{3}+e_{3}\otimes e_{1};\ \varepsilon _{3,3}^{3}\left( e_{1}\right)
=1;\ \varepsilon _{3,3}^{3}\left( e_{2}\right) =1;\ \varepsilon
_{3,3}^{3}\left( e_{3}\right) =0. $
\end{enumerate}

For the multiplication $\mu _{4}^{3}$,  there does not exist any
bialgebras.

For the multiplication $\mu _{5}^{3}$, we have

\begin{enumerate}
\item  $ \Delta _{5,1}^{3}\left( e_{1}\right) =e_{1}\otimes
e_{1};\ \Delta _{5,1}^{3}\left( e_{2}\right) =e_{2}\otimes
e_{2};\ \Delta _{5,1}^{3}\left( e_{3}\right) =e_{2}\otimes
e_{3}+e_{3}\otimes e_{2};\ \varepsilon _{5,1}^{3}\left( e_{1}\right)
=1;\ \varepsilon _{5,1}^{3}\left( e_{2}\right) =1;\ \varepsilon
_{5,1}^{3}\left( e_{3}\right) =0. $
\end{enumerate}

In the sequel we consider that all the algebras are unital and the unit $\eta$  corresponds to $e_1$.
\begin{prop}
The bialgebras on $\K ^3$ which are unital infinitesimal bialgebras are given by the following pairs of multiplication and comultiplication with the appropriate unit and counits.
$$\left( \mu _{1}^{3},\Delta _{1,2}^{3}\right) ,\ \left( \mu
_{1}^{3},\Delta _{1,5}^{3}\right) ,\ \left( \mu _{1}^{3},\Delta
_{1,6}^{3}\right) ,\ \left( \mu _{1}^{3},\Delta _{1,8}^{3}\right)
,\ \left( \mu _{1}^{3},\Delta _{1,11}^{3}\right) ,\ \left( \mu
_{1}^{3},\Delta _{1,14}^{3}\right),\ \left(\mu _{1}^{3},\Delta _{1,15}^{3}\right) ,\ $$
$$\left( \mu
_{1}^{3},\Delta _{1,18}^{3}\right) ,\ \left( \mu _{2}^{3},\Delta
_{2,1}^{3}\right) ,\ \left( \mu _{2}^{3},\Delta _{2,2}^{3}\right) ,\ \left( \mu
_{3}^{3},\Delta _{3,2}^{3}\right) ,\ \left( \mu _{3}^{3},\Delta
_{3,3}^{3}\right) ,\ \left(\mu _{5}^{3},\Delta _{5,1}^{3}\right). $$
\end{prop}

\begin{rem}
The unital infinitesimal bialgebras found here could be used to produce examples of Rota-Baxter algebras.
\end{rem}

We summarize, in the following table, the numbers of non-isomorphic
bialgebras and unital infinitesimal bialgebras associated to a given
algebra.

\vspace{0.2cm}
\begin{center}
 \begin{tabular}{|c|c|c|}
   \hline
    \textit{algebra}& \textit{bialgebras} & \textit{infinitesimal bialgebras} \\ \hline
   $\mu _{1}^{3}$ & 18 & 8 \\ \hline
   $\mu _{2}^{3}$ & 3 & 2 \\ \hline
   $\mu _{3}^{3}$ & 3 & 2 \\ \hline
   $\mu _{4}^{3}$ & 0 & 0 \\ \hline
   $\mu _{5}^{3}$ & 1 & 1 \\
   \hline
 \end{tabular}
 \end{center}
 \vspace{0.2cm}

\begin{cor}
\

\begin{itemize}
\item The number of 3-dimensional trivial 2-associative bialgebras
 (same multiplication) is \textbf{13}.

\item We have the following 3-dimensional non-trivial 2-associative
bialgebras   (different multiplication)

 $\left( \mathbb{K}^3,\mu _{3}^{3},\mu _{5}^{3},\eta,\Delta
_{3,1}^{3},\varepsilon_{3,1}^{3}\right) ,\left( \mathbb{K}^3,\mu _{1}^{3},\mu
_{2}^{3},\eta,\Delta _{2,1}^{3},\varepsilon_{2,1}^{3}\right) $ and $\left( \mathbb{K}^3,\mu
_{1}^{3},\mu _{2}^{3},\eta,\Delta _{2,2}^{3},\varepsilon_{2,2}^{3}\right) .$
\end{itemize}
\end{cor}

 \begin{prop}
 The number of 3-dimensional isomorphism classes of 2-bialgebras is
\begin{itemize}
\item  type (1,1) : 25 .

\item type (1,2) : 159.

\item type (2,1) : 1. Namely  $\left( \mathbb{K}^3,\mu _{3}^{3},\mu
_{5}^{3},\eta,\Delta _{3,1}^{3},\Delta _{5,1}^{3},\varepsilon
_{3,1}^{3},\varepsilon _{5,1}^{3}\right) $.

 \item  type (2,2) : 3. Namely,
$\left( \mathbb{K}^3,\mu _{1}^{3},\mu _{2}^{3},\eta,\Delta
_{1,3}^{3},\Delta _{2,1}^{3},\varepsilon _{1,3}^{3},\varepsilon _{2,1}^{3}\right)$ ,

$\left( \mathbb{K}^3,\mu _{1}^{3},\mu _{2}^{3},\eta,\Delta
_{1,4}^{3},\Delta _{2,1}^{3},\varepsilon _{1,4}^{3},\varepsilon _{2,1}^{3}\right)$ and
$\left( \mathbb{K}^3,\mu _{1}^{3},\mu _{2}^{3},\eta,\Delta
_{1,5}^{3},\Delta _{2,1}^{3},\varepsilon _{1,5}^{3},\varepsilon _{2,1}^{3}\right). $
\end{itemize}
\end{prop}
\begin{cor}
There exists only one 2-2-bialgebra of dimension 3 which is given by
$$\left( \mathbb{K}^3,\mu _{1}^{3},\mu _{2}^{3},\eta,\Delta
_{1,5}^{3},\Delta _{2,1}^{3},\varepsilon _{1,5}^{3},\varepsilon _{2,1}^{3}\right). $$
\end{cor}

\end{document}